\documentclass[12pt]{article}

\usepackage{times}
\usepackage{arydshln}

\usepackage{amsmath,amsthm,amssymb}
\usepackage{graphicx}
\usepackage{color}
\usepackage{url}
\usepackage{multirow}

\def \bx {\pmb{x}}
\def \bvf {\pmb{f}}
\def \bh {\pmb{h}}
\def \bw {\pmb{w}}
\def \bs {\pmb{s}}

\def \deta {\pmb{\eta}}

\definecolor{OliveGreen}{rgb}{0.3, 0.6, 0.1}

\title{One-way dependent clusters and stability of cluster synchronization in directed networks} 

\author
{Matteo Lodi,$^{1}$ Francesco Sorrentino,$^{2}$ Marco Storace $^{1\ast}$\\
	\\
	\normalsize{$^{1}$DITEN, University of Genoa, Via Opera Pia 11a, I-16145, Genova, Italy}\\
	\normalsize{$^{2}$Mechanical Engineering Department, University of New Mexico,}\\
	\normalsize{Albuquerque, NM 87131, USA}\\
	\\
	\normalsize{$^\ast$To whom correspondence should be addressed; E-mail:  marco.storace@unige.it.}
}

\date{}

\begin{document} 
	
	\baselineskip24pt
	
	\maketitle

\begin{abstract}
	\textbf{Cluster synchronization in networks of coupled oscillators is the subject of broad interest from the scientific community, with applications ranging from neural to social and animal networks and technological systems. Most of these networks are directed, with flows of information or energy that propagate unidirectionally from given nodes to other nodes. Nevertheless, most of the work on cluster synchronization has focused on undirected networks. Here we characterize cluster synchronization in general directed networks. Our first observation is that, in directed networks, a cluster A of nodes might be one-way dependent on another cluster B: in this case, A may remain synchronized provided that B is stable, but the opposite does not hold. The main contribution of this paper is a method to transform the cluster stability problem in an irreducible form. In this way, we decompose the original problem into subproblems of the lowest dimension, which allows us to immediately detect inter-dependencies among clusters. We apply our analysis to two examples of interest, a human network of violin players executing a musical piece for which directed interactions may be either activated or deactivated by the musicians, and a multilayer neural network with directed layer-to-layer connections.}
\end{abstract}

\vspace{2cm}

\section*{Introduction}\label{sec:intro}
Synchronization is a pervasive phenomenon in networks of natural and man-made dynamical systems \cite{pikovsky:2003}, which often enables complex functions corresponding to patterns of coordinated activity. Some of these systems require complete (or full) synchronization among all network components
to properly function: examples include the coherent firing of neural populations \cite{shlens:2008}, synchronous behaviors in networks of chemical oscillators \cite{kiss:2002}, the concerted rhythmical flashing on and off of firefly swarms \cite{stone:2002}, or frequency-locked power generation \cite{grainger:2003}. In other cases, such as higher brain centers of the nervous system \cite{guevara:2017}, rhythmic activity and synchronization may correspond to a pathological condition. Other systems rely on cluster (or partial) synchronization (CS), where different groups exhibit different, yet synchronized, internal behaviors. For instance, a large corpus of studies has shown the connection of neural synchronization of brain areas with cognition \cite{guevara:2017}, animal locomotion is controlled to some extent by central pattern generators that impose a given gait by exploiting synchronous clusters of neurons \cite{grillner:2006,Ijspeert:2008,goulding:2009}, and connected vehicle systems often give rise to CS \cite{orosz:2009}.
The study of CS is therefore relevant to analyze and control natural (ecological, immune, neural, and cellular) \cite{kaneko:1994} and artificial \cite{buono:2015,buono:2018,kiss:2019} systems.

A network of dynamical systems can be modeled through (i) a graph (where an element is called a node and a connection is called an edge or link, or an arrow if oriented), which is an abstract description of its architecture or topology, and (ii) a set of dynamical equations describing the time evolution of each node and, possibly, also of each connection. 
The interplay between the network topology and the dynamics of the nodes and of their interactions determines the complex dynamics emergent in all the systems described above.

A large body of literature has investigated stability of the cluster synchronous solution for \textit{undirected networks} \cite{golubitsky:2005,golubitsky:2006,golubitsky:2012,pecora:2014,aguiar:2014,schaub:2016,sorrentino:2016,siddique:2018,dellarossa:2020,zhang:2020}, corresponding to graphs whose edges do not have a direction and indicate a two-way relationship. However, most real networks are \textit{directed} and are described by graphs (also called digraphs) whose edges have a direction that indicates a one-way relationship. In many fields, from systems biology to the engineering of distributed technology \cite{kirst:2016,hens:2019}, the corresponding physical systems do not always exhibit a bidirectional flow of information/energy, thus leading to preferred paths for routing of information/energy. For example, neural, gene regulatory, metabolic, epidemic networks are all directed. Understanding the conditions for the emergence of stable synchronized clusters in these networks is of considerable importance.

Following Ref. \cite{pecora:2014} a fundamental issue is to be able to reduce the dimensionality of the stability analysis for synchronized clusters. However, it remains an open problem how such reduction can be obtained for directed networks. Here we bridge this gap by decoupling the CS stability problem into lower-dimensional problems. References \cite{pecora:2014,sorrentino:2016,siddique:2018} study undirected networks and explain how to analyze the stability of any CS solution by introducing a transformation - based on the \textit{irreducible representation (IRR)} of the network symmetry group - from the node coordinate system to the IRR coordinate system.

The methods based on IRR transformations  \cite{golubitsky:2012,pecora:2014,dellarossa:2020} have been applied to undirected networks. For these networks, it is also possible to understand if two clusters $C_1$ and $C_2$ are \textit{intertwined}, i.e., if desynchronization of cluster $C_1$ leads to loss of synchrony in cluster $C_2$ and \textit{vice versa} \cite{pecora:2014,cho:2017}.

The main contribution of this paper is a mathematical method for analyzing the stability of synchronized clusters in directed multilayer networks with general structure. This method is based on two cornerstones. The first one is a generalization of the concept of intertwined clusters, which takes into account directionality: we define \textit{one-way dependent clusters}, meaning that the breaking (i.e., loss of stability) of a synchronized cluster $C_1$ causes the loss of synchrony of another cluster $C_2$, but not \textit{vice versa}. The second cornerstone is a method to build up the matrices that we use to analyze cluster stability in directed networks through linearization: owing to this method, we solve a problem still open, namely how to find a coordinate transformation that separates \textit{as much as possible} perturbation modes in the stability analysis.
This method also allows easy detection (by inspection of certain matrices) of inter-dependencies between clusters in directed networks, and so a classification of pairs of clusters as independent or intertwined or one-way dependent in terms of their stability. This method enables the analysis of a class of systems, for which, to the best of our knowledge, other analysis tools are not available, thus bridging a gap towards the analysis of biological, social and technological systems described by directed networks. Alternative changes of variables \cite{cho:2017} are simpler to compute but do not ensure that the coordinate transformation is optimal in the sense specified above.

We apply the proposed method to two illustrative case studies, namely two networks from different fields. The first one is a single-layer network with delays and delay-dependent connections modeling a group of violin players \cite{shahal:2020}. The second one is a two-layer network of neurons that can also produce chimera states \cite{majhi:2017,ruzzene:2020}. In both cases our method provides a key to better understand observed phenomena in terms of cluster synchronization.

\section*{Results}\label{sec:results}

A network can be represented by a weighted digraph (or directed graph), defined by a vertex set $\mathcal{V} = \{v_1, \ldots, v_N\}$ (representing the nodes), an arrow set $\mathcal{E}_k$ for each kind of link ($k= 1,\ldots,L$), and a weight set $\mathcal{W}_k$ representing the connection weights for the $k$-th kind of link. The nodes can be of $M$ different kinds.

The \textit{weighted} adjacency matrix $A^k$ embeds information about the arrow set $\mathcal{E}_k$ and the weight set $\mathcal{W}_k$: $A_{ij}^k \in \mathcal{W}_k$ is the (nonzero) weight of the link going from $v_j$ to $v_i$; $A_{ij}^k = 0$ if there is no arrow from $v_j$ to $v_i$ and $A_{ij}^k$ belongs to $\mathcal{W}_k$ otherwise.
In an \textit{undirected} network, if there is an arrow with a certain weight going from $v_i$ to $v_j$ then there is also an arrow with the same weight going from $v_j$ to $v_i$; in this case, $A^k$ is a symmetric matrix and the underlying graph is undirected.

A \textit{multilayer} network is characterized by the presence of different types of nodes interacting through different types of connections, with each \textit{layer} being formed of nodes of the same type and/or connections of the same type \cite{boccaletti:2014,kivela:2014}. In particular, in this paper a layer is given by nodes of the same type \cite{dellarossa:2020}.
Arrows can be inside each layer (intralayer connections) or go from layer to layer (interlayer connections).

\begin{figure}[h!]
	\centering
	\includegraphics[width=1\textwidth]{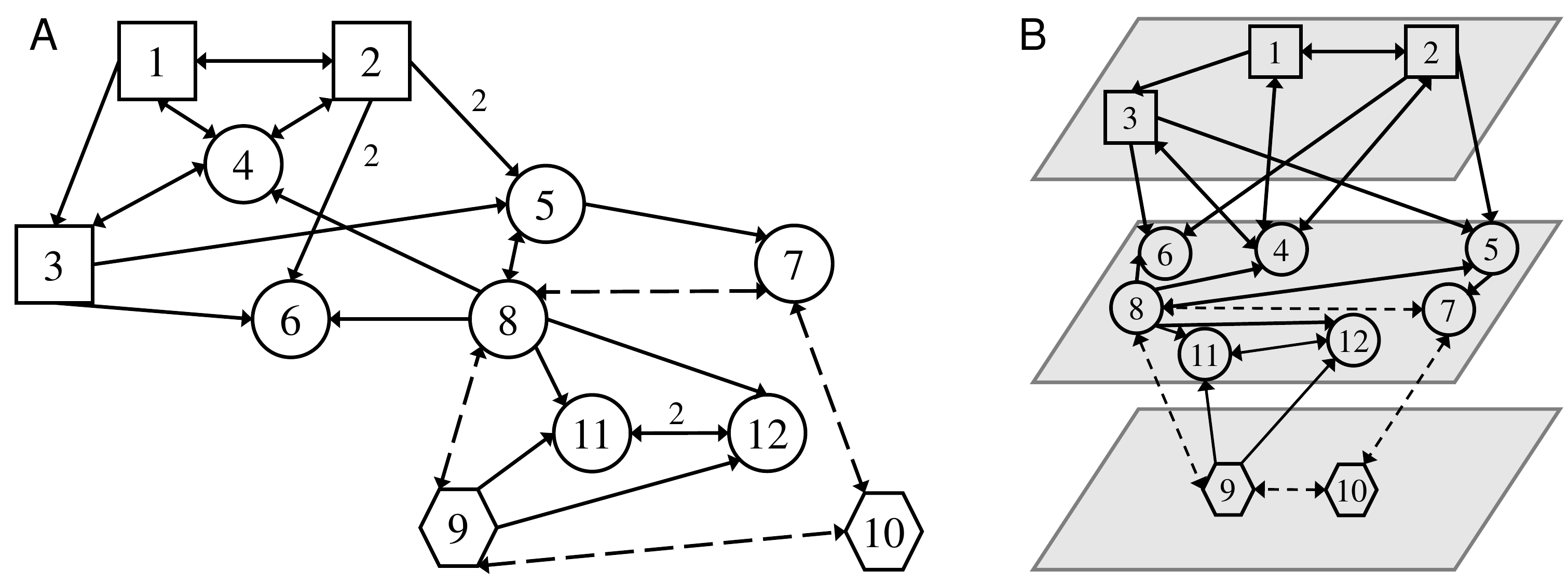}
	\caption{\footnotesize{Directed multilayer network. (A) Directed network with $N=12$ nodes of $M=3$ shape-coded different kinds and $L=2$ kinds of links, represented by either solid or dashed arrows. The vertex set is $\mathcal{V} = \{1, \ldots, 12\}$, the arrow sets are $\mathcal{E}_1= \{ (1,2), (1,3), (1,4), (2,1), (2,4), (2,5), (2,6), \ldots\}$ (solid arrows) and $\mathcal{E}_2= \{ (7,8), (7,10), (8,7), (8,9), \ldots \}$ (dashed arrows). All weights in the sets $\mathcal{W}_1$ and $\mathcal{W}_2$ are equal to 1 except for the edges (2,5), (2,6), (11,12), and (12,11) (all belonging to $\mathcal{E}_1$) that have weight equal to 2. (B) The three layers for the network in A: top layer with square nodes, middle layer with circle nodes, and bottom layer with hexagon nodes. Examples of intralayer connections: $1 \rightarrow 2$ (top layer), $8 \leftrightarrow 5$ (middle layer), $9 \leftrightarrow 10$ (bottom layer). Examples of interlayer connections: $2 \rightarrow 5$ (between top and middle layers), $8 \leftrightarrow 9$ (between middle and bottom layers).}}
	\label{fig:exnetw}
\end{figure}

The above definitions are illustrated through the network shown in Fig. \ref{fig:exnetw}.

We now associate a dynamics to each node, by considering a network with $N$ nodes, which are connected through $L$ different kinds of links. We model networks of this kind through the following set of dynamical equations ($i = 1,\ldots,N$) \cite{lodi:2020}, which provides a rather general description of deterministic systems governed by pairwise interactions
\begin{equation}
\dot \bx_i =  \bvf_i(\bx_i(t)) + \sum_{k=1}^{L}\sigma_k \sum_{j=1}^{N} A_{ij}^k \bh^k(\bx_i(t),\bx_j(t-\delta_k)),
\label{eq:dyneq}
\end{equation}
where $\bx_i \in \mathbb{R}^n$ is the $n$-dimensional state vector of the $i$-th node, $\bvf_i : \mathbb{R}^n \rightarrow \mathbb{R}^n$ is the vector field of the isolated $i$-th node describing the uncoupled dynamics of $\bx_i$, $\sigma_k \in \mathbb{R}$ is the coupling strength of the $k$-th kind of link, $A^k$ is the weighted adjacency matrix with respect to the $k$-th kind of link, for which the pairwise interaction between two generic nodes $i$ and $j$ is described by the nonlinear function $\bh^k: \mathbb{R}^n \times \mathbb{R}^n \rightarrow \mathbb{R}^n$, and $\delta_k$ is the transmission delay characteristic of the $k$-th kind of link.
We assume each individual node can be of one out of $M$ different types (with $M \leq N$): $\bvf_i(x) = \bvf_j(x)$ if $i$ and $j$ are of the same type, $\bvf_i(x) \neq \bvf_j(x)$ otherwise. Within this general framework, where all oscillators can be different, if $M << N$ the vector fields $\bvf_i$ are not all different, but belong to a restricted set of $M$ models.
Equation \eqref{eq:dyneq} describes the dynamics of a directed multilayer network with both nodes of different types and arrows of different types. The interpretation of Eq. \eqref{eq:dyneq} in terms of layers is made explicit in the \textit{Supplementary Material}, Sec. 1.

This paper deals with cluster synchronization of the system \eqref{eq:dyneq}, i.e., ${\bx}_i(t) = \bx_j(t)$ for $i$ and $j$ that belong to the same cluster.

\subsection*{Existence}
The \textit{existence} of a synchronized cluster can be either ruled out or not, based on the analysis of the \textit{network topology} and of the isolated node vector fields $\bvf_i$.
By a \textit{partition} $\mathcal{P} = \{C_1, C_2,\ldots,C_Q\}$ of a graph, we mean a partition of its vertex set into clusters which satisfies the following properties: $\mathcal{P} = \left\{ C_i \subset \mathcal{V} : C_i \cap C_j = \emptyset \quad \forall i\neq j, \quad \displaystyle{\bigcup_{i=1}^{Q}} C_i = \mathcal{V} \right\}$, where we call $C_i$ the $i$-th \textit{cluster}, with $i = 1,\ldots,Q$. We can identify each cluster through the labels of its constituting vertices and a given color.
A cluster is said to be \textit{trivial} if it consists of only one node.

Let $N_i$ be the number of nodes belonging to the cluster $C_i$. Therefore, $\sum_{i=1}^{Q} N_i = N$.

Each graph can admit different partitions, corresponding to different colorings of the graph.
\textit{Invariant partitions}, also called \textit{balanced partitions} or \textit{equitable partitions}, have a large number of applications in many fields such as graph theory, network theory and control theory (see \cite{neuberger:2020} and references therein). 
An invariant partition is one in which nodes with the same color follow the same time evolutions. Consequently, a partition (and the corresponding coloring) is balanced/equitable if all nodes with label (color) $p$ have the same vector field $\bvf_i$ and get the same overall input with delay $\delta_k$ from the nodes of label (color) $q$, for $p,q = 1,\ldots,Q$ and for any delay $\delta_k$ with $k = 1,\ldots,L$, namely if
\begin{equation} \label{eq:equitablecluster}
\left\{ \begin{array}{c} \displaystyle{\sum_{v_a \in C_q} A_{ia}^k = \sum_{v_a \in C_q} A_{ja}^k}\\ \bvf_i = \bvf_j \end{array}  , \quad \begin{array}{l} \forall \; v_i,v_j \in C_p\\ \forall \; C_p, C_q \in \mathcal{P} \end{array} \right. .
\end{equation}

The existence of a balanced/equitable partition (coloring) for a given network is a necessary condition for the existence of a CS solution \cite{golubitsky:2006}.

In particular, among all possible equitable partitions, we focus on the partition $\mathcal{P}$ that corresponds to the so-called \textit{minimal balanced coloring}, namely a balanced coloring with the minimal number of colors \cite{mckay:1981,belykh:2011,lodi:2020ISCAS}. We also remark that the minimal balanced coloring provides the cluster configuration with the minimum number $Q$ of clusters.

\begin{figure}[h!]
\includegraphics[width=1\textwidth]{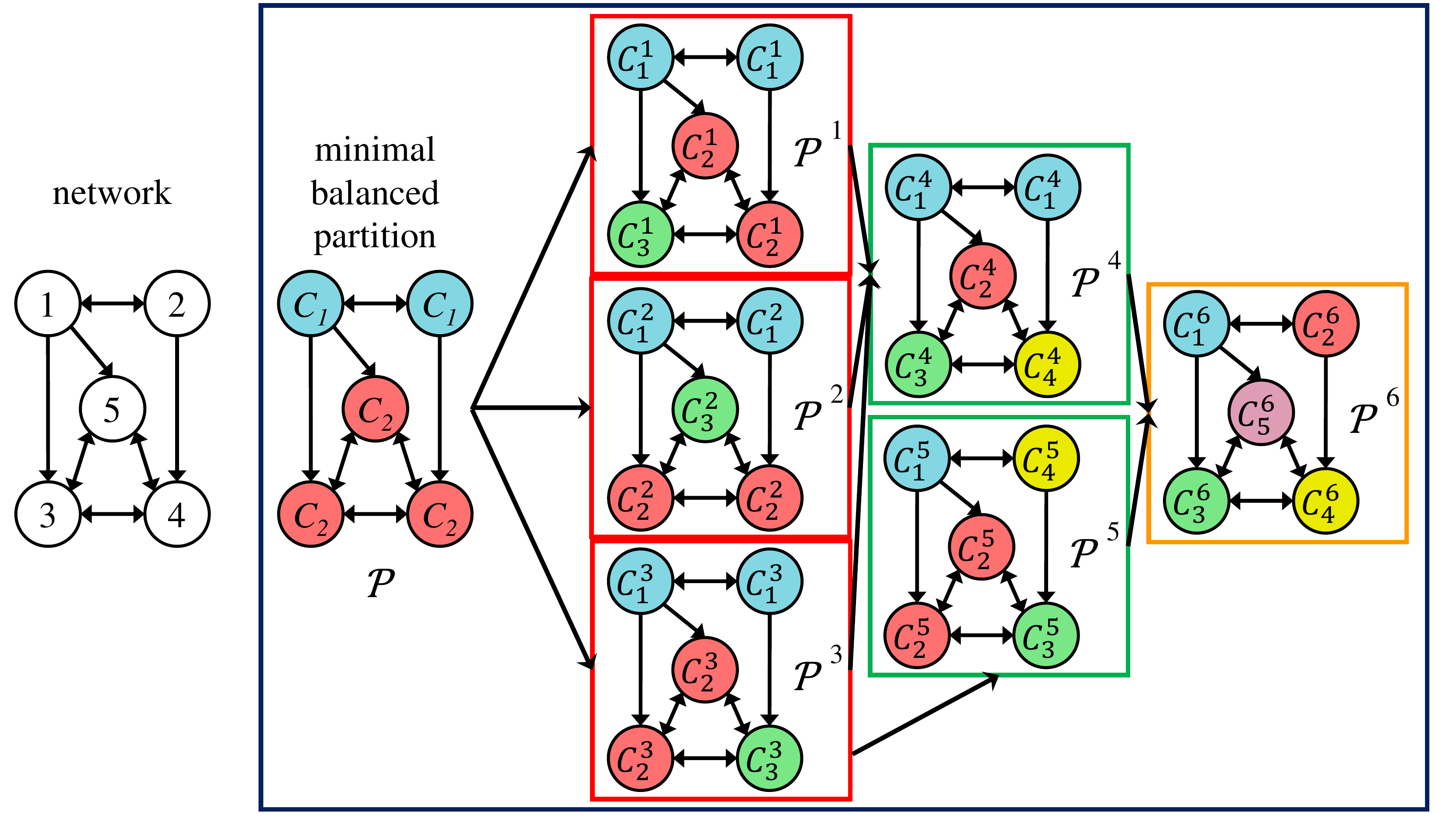}
	\caption{\footnotesize{Illustrative example of balanced partitions. Balanced partitions (within the blue rectangle) of a simple directed network with $N=5$ nodes. The partition $\mathcal{P}$ corresponds to the minimal balanced coloring, with two clusters $C_1$ (containing nodes 1 and 2, colored blue) and $C_2$ (containing nodes 3, 4, and 5, colored red): each node of cluster $C_2$ receives one arrow from one node of $C_1$; each node of cluster $C_1$ receives one arrow from one node of $C_1$ and two arrows from nodes of $C_2$. Breaking these clusters, we can obtain all the other balanced partitions, which are non-minimal, because they contain a higher number of clusters. Each partition $\mathcal{P}^1$,  $\mathcal{P}^2$, and $\mathcal{P}^3$ (within red squares) contains three clusters, which are obtained from $\mathcal{P}$ by breaking the largest cluster into two smaller ones. For instance, $\mathcal{P}^1$ contains $C_1^1$ (with nodes 1 and 2, colored blue), $C_2^1$ (with nodes 4 and 5, colored red), and $C_3^1$ (with node 3, colored green). This partition is balanced, because each node of $C_1^1$ receives one arrow from one blue node, each node of $C_2^1$ receives one arrow from a blue node, one arrow from a red node and one from a green node, and each node of $C_3^1$ receives one arrow from a blue node and two arrows from red nodes. Similarly, partitions $\mathcal{P}^4$ and $\mathcal{P}^5$ (within green squares) contain four clusters each and are obtained from the partitions with three clusters by breaking the largest cluster into two smaller clusters: the partition $\mathcal{P}^4$ is obtained by breaking the red clusters of $\mathcal{P}^1$,  $\mathcal{P}^2$, or $\mathcal{P}^3$; the partition $\mathcal{P}^5$ is obtained by breaking the blue cluster of $\mathcal{P}^3$. Finally, the partition $\mathcal{P}^6$ (within orange square) contains five trivial clusters (each one contains only one node).}}
	\label{fig:balpart}
\end{figure}

Figure \ref{fig:balpart} shows the seven balanced partitions of a simple single-layer directed network, whose adjacency matrix is reported in the \textit{Supplementary Material}, Sec. 1. The partition labeled $\mathcal{P}$ corresponds to the minimal balanced coloring. The other partitions (labeled $\mathcal{P}^j$) correspond to non-minimal balanced colorings and are composed of clusters $C_i^j$. In particular, $\mathcal{P}^6$ contains only trivial clusters. It is easy to check that the partitions in Fig. \ref{fig:balpart} are all of the balanced partitions for this network.

Each balanced partition admits a cluster synchronous solution where all the nodes in the same cluster follow the same time evolution. This solution can be expressed in terms of the so-called \textit{invariant synchrony subspaces} (ISS), i.e., the polydiagonal subspaces under all admissible vector fields \cite{stewart:2003,aguiar:2014,steur:2016,aguiar:2017,aguiar:2018,neuberger:2020}.

The ISS $\Delta^j$ corresponding to each partition $\mathcal{P}^j$ of Fig. \ref{fig:balpart} can be easily determined by inspection: $\Delta^1$ is such that the state evolution of nodes 4 and 5 is the same, i.e., $\{\bx_4=\bx_5\}$; $\Delta^2$ is such that $\{\bx_1=\bx_2 \; \& \; \bx_3=\bx_4\}$; $\Delta^6$ is such that each node of the network has a different time evolution, and so forth.

We remark that not all the balanced partitions compatible with a given topology can also be observed when taking into account the network dynamics \eqref{eq:dyneq}, as some of these may be unstable. Thus the dynamics governing the network must be incorporated in the study of stability for the CS solution.

\subsection*{Stability}
The possibility of observing a synchronized cluster in a given physical system depends on its stability, and hence on the particular dynamics of the nodes and connections.

Once the $Q$ equitable clusters corresponding to a balanced partition are found, the CS analysis can be focused on a simplified dynamical model (called \textit{quotient network}) whose $Q$ nodes correspond to each one of the clusters. In particular, we analyze the cluster stability by linearizing Eq. \eqref{eq:dyneq} about a state corresponding to exact synchronization among all the nodes within each cluster \cite{lodi:2020}.
The Master Stability Function (MSF) approach \cite{pecora1998master}, which evaluates the stability of synchronized solutions through linearization, is based on two standard steps: (i) finding the variational equations of the network about the synchronized solutions and (ii) expressing these variational equations in a new system of coordinates, which decouples the perturbation dynamics along the transverse manifold from that along the synchronous manifold.

We collect all state trajectories in the vector ${\bx}(t) = [\bx_1^T(t), \bx_2^T(t),\dots,\bx_N^T(t)]^T$. As it is possible for all the nodes within a cluster to synchronize, we define the $q$-th cluster state: $\bs_q(t) = \bx_i(t)$ for all $i$ in cluster $C_q$. Correspondingly, the network can produce $Q$ distinct synchronized motions $\{\bs_1(t), \bs_2(t),\ldots,\bs_Q(t)\}$, one per cluster. We collect them in the vector ${\bs}(t) = [\bs_1^T(t), \bs_2^T(t),\dots,\bs_Q^T(t)]^T$.
We then analyze the dynamics of a small perturbation $\bw_i(t) = \bx_i(t)-\bs_{q_i}(t)$ ($i=1,\ldots,N$), where $\bs_{q_i}(t)$ is the cluster state for node $i$ in cluster $C_q$, by linearizing around a specific network solution ${\bs}(t)$.

In order to reduce the dimensionality of the stability analysis, we first divide the full state space of the variational equation into minimal flow-invariant subspaces and then calculate the maximum Lyapunov exponent (MLE) in each subspace, in order to determine whether perturbations within that subspace grow or decay. A key aspect for an efficient analysis of the stability of arbitrary synchronization patterns \cite{irving:2012,pecora:2014,zhang:2020} is to be able to optimally decouple the perturbation modes.
To this end, we make a change of coordinates through the transformation matrix $T$, which maps the perturbation vector $\bw = [\bw_1^T,\bw_2^T,\ldots,\bw_N^T]^T$ into a new coordinate system $\deta = [\deta_1^T,\deta_2^T,\ldots,\deta_N^T]^T = (T \otimes \mathbb{I}_n) \bw$, where the symbol $\otimes$ denotes the tensor product of matrices, i.e., the Kronecker product.
The $N \times N$ matrix $T$ is made up of two parts: $T = \begin{bmatrix} T_\parallel \\ T_\perp \end{bmatrix}$, where the $Q \times N$ submatrix $T_\parallel$ is concerned with the directions inside the synchronization manifold and the corresponding perturbations do not influence the stability of the synchronized clusters; on the other hand, the $(N-Q) \times N$ submatrix $T_\perp$ is concerned with the directions transverse to the synchronization manifold and the evolution of the variational equation along these directions determines the stability of the synchronized clusters. Therefore, the perturbation dynamics are separated into two parts: that along the synchronous manifold (described by the first $Q$ components $\deta_i$), which accounts for the stationary/periodic/chaotic state dynamics within each cluster and that transverse to it (described by the last components $\deta_i$, $i\in [Q+1, N]$), which accounts for the stability of clusters.

In this perspective, the coordinate system $\deta$ must separate the evolution of transverse and parallel perturbation modes, by decoupling the transverse perturbation modes as much as possible. In particular, the entries of the $q$-th row of $T_\parallel$ are all zero, except for those whose column index denotes one of the $N_q$ nodes belonging to $C_q$. The value of the non-zero entries is $1/\sqrt{N_q}$.

On the other hand, $T_\perp$ has not a unique structure but must satisfy two constraints:
\begin{itemize}
	\item (A) $T_\perp$ contains $N_q - 1$ rows for each cluster $C_q$; the $N$ entries of each row are are all zero, with the exception of those whose column index denotes one of the $N_q$ nodes belonging to $C_q$;
	\item (B) each row of $T_\perp$ is orthonormal to any other row of $T$.
\end{itemize}

The matrix $T$ transforms the coupling matrices $A^k$ into matrices $B^k=T A^k T^{-1}$. 
The system of equations governing the dynamics of the perturbation vector $\bw$ is \textit{reducible}, because all the adjacency matrices $A^k$ can be put in block-upper triangular form $B^k = T A^k T^{-1}$ by application of the same invertible transformation matrix $T$, where each matrix $B^k$ has the block structure\\
\begin{equation}
	B^k = \begin{bmatrix} B^k_{11} & B^k_{12}\\ 0 & B^k_{22}
	\end{bmatrix}
\end{equation}
for any $k$. Note that this corresponds to effectively decoupling the perturbation vector into two parts, where the second one (containing the transverse perturbation modes) evolves independently of the first one (containing the perturbations inside the synchronization manifold). Among the many possible matrices $T$, we look for the one that ensures the maximal separation of the transverse perturbation modes.

A coordinate system $\deta$ is no longer reducible (i.e., it is \textit{irreducible}), if it is impossible to further separate the transverse perturbation modes $\deta_i$ by application of a different transformation matrix $T$. We then call irreducible the corresponding set of matrices $B^k$ that cannot be reduced further.

In the case of undirected networks, we can find a transformation matrix $T$ such that $B^k_{12} = 0$ (null block) for any $k$. Therefore, each matrix $B^k$ can be reported in block-diagonal form and there are methods -- based on either the group theory \cite{tinkham:2003,siddique:2018,dellarossa:2020} or matrix algebras \cite{irving:2012,zhang:2020} -- that transform it in an irreducible form.

In the case of directed networks, the matrix $B^k_{12}$ is not a null block, in general. Therefore, each matrix $B^k$ can be reported in block-upper triangular form \cite{cho:2017,lodi:2020}, but, at the best of our knowledge, there are no methods to find the matrix $T$ that transforms the perturbation modes in an irreducible form. We remark that the matrix must be decomposed so as to separate perturbations along the synchronous manifold from those along the transverse manifold. This prevents the use of standard matrix decomposition techniques, which do not take this into account.

We aim to find the matrix $T$ that converts the node coordinate system into the irreducible coordinate system, thus evidencing the inter-dependencies among the perturbation components in synchronized clusters of directed networks. Our proposal of a method to build up this matrix for directed networks is the major contribution of this paper.

\begin{figure}[b!!]
	\centering
	\includegraphics[width = 0.9\linewidth]{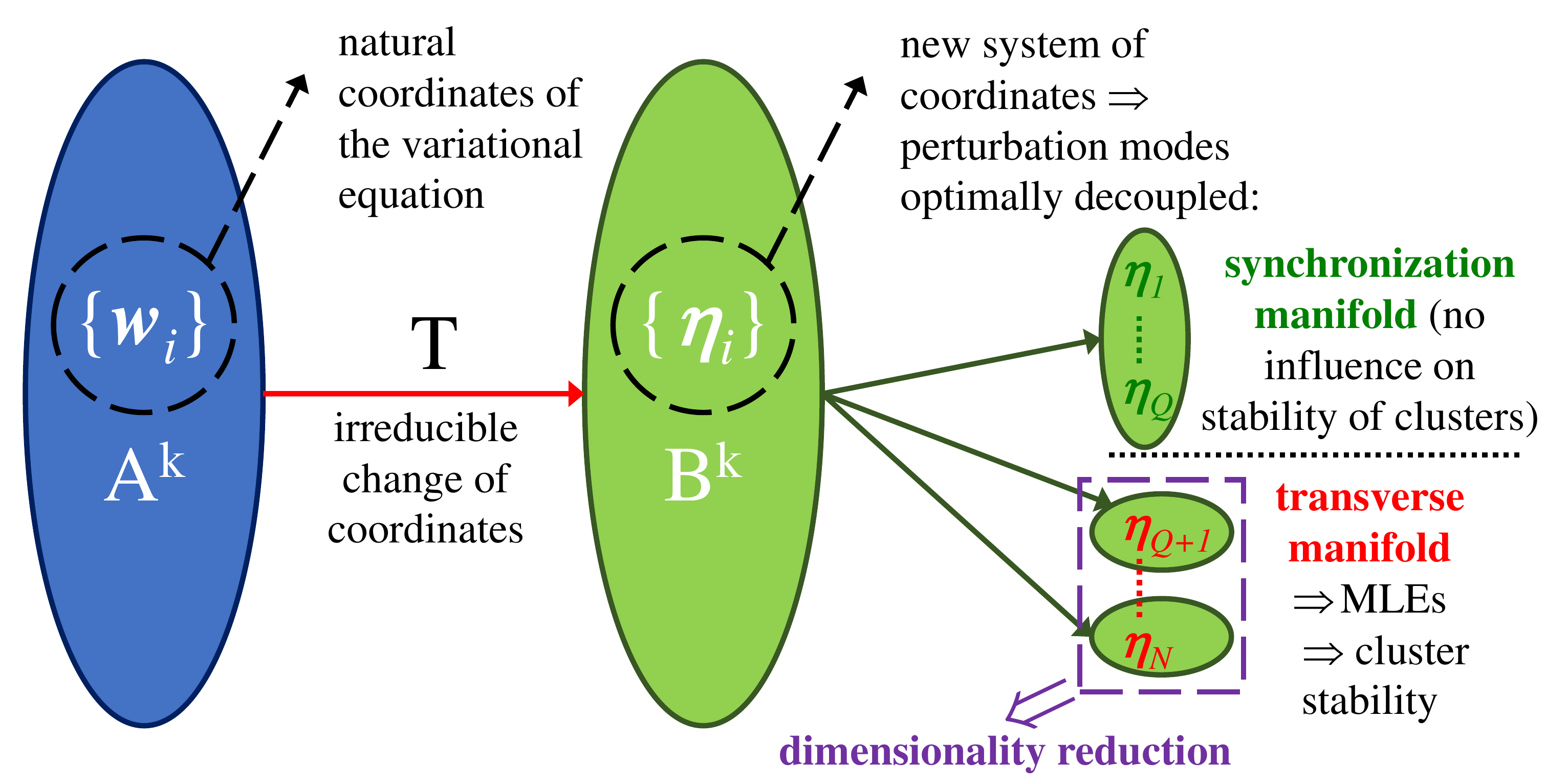}
	\caption{\footnotesize{Dimensionality reduction based on the transformation matrix $T$. First (blue oval), the system is linearized about a synchronous solution $\bs(t)$ (composed of $Q$ distinct synchronized motions), in order to analyze the dynamics of a small perturbation $\bw_i(t)$ around each node $i$ ($i = 1, \ldots, N$); the topology of each network layer is described by an adjacency matrix $A^k$. Second (large green oval), the transformed system is described by a new set of coordinates $\{ \deta_i \}$ which are coupled through the matrices $B^k = T A^k T^{-1}$. The new coordinates are decoupled into two sets (small green ovals), one containing the perturbations lying on the synchronization manifold, which do not influence the stability of the synchronized clusters, and one containing the perturbations transverse to the synchronization manifold, which influence the stability of the synchronized clusters. In particular, we look for the irreducible change of coordinates (and corresponding matrix $T$) which decouples the transverse perturbation modes as much as possible (green ovals within the dashed purple rectangle). The transverse perturbation modes contain information on the cluster stability, also evidencing the presence of inter-dependencies between different clusters.}}
	\label{fig:IRRdecoupling}
\end{figure}

The process for optimally decoupling the perturbation modes described above is summarized in Fig. \ref{fig:IRRdecoupling}.

\subsection*{One-way dependent clusters and breaking vectors}\label{sec:dirintertw}
As stated in the Introduction, in undirected networks two clusters $C_1$ and $C_2$ are said to be \textit{intertwined} \cite{pecora:2014,cho:2017} if  desynchronization of cluster $C_1$ can lead to loss of synchrony in cluster $C_2$ and \textit{vice versa}. In directed networks, instead, it may happen that the stability of $C_1$ depends on the stability of $C_2$ but the opposite is not necessarily true. In this case, we say that $C_1$ is \textit{one-way dependent} on $C_2$.

For instance, the equitable partition $\mathcal{P}$ in Fig. \ref{fig:balpart} contains two clusters, $C_1 = \{1,2\}$ and $C_2 = \{3,4,5\}$. The stability of $C_1$ influences also the stability of $C_2$: if only $C_1$ desynchronizes, then its nodes have no longer the same color and therefore also the nodes of $C_2$ cannot have the same color; in other words, they are no longer synchronized, which corresponds to either partition/coloring $\mathcal{P}^1$ or $\mathcal{P}^2$. On the other hand, if only $C_2$ desynchronizes, in all the possible corresponding balanced partitions ($\mathcal{P}^3$, $\mathcal{P}^4$, $\mathcal{P}^5$, $\mathcal{P}^6$) the nodes of $C_1$ keep the same color.
Therefore $C_2$ is one-way dependent on $C_1$.

Our strategy to check for the presence of one-way dependent clusters -- also when the network has a large number of nodes and one cannot easily detect the inter-dependencies between clusters by inspection -- is based on the ISS. We consider an equitable cluster belonging to the partition $\mathcal{P}$: any other balanced partition $\mathcal{P}^j$ can be obtained by breaking one or more clusters of $\mathcal{P}$ into the subclusters $C_i^j$.

We consider the minimal balanced partition $\mathcal{P}$, the non-minimal balanced partition $\mathcal{P}^j$, and the $i$-th cluster $C_i^j \in \mathcal{P}^j$.
Let $\mathcal{B}_q^j$ be the $N$-length row vector describing the \textit{transition} from the cluster $C_q \in \mathcal{P}$ to smaller clusters $C_i^j \in \mathcal{P}^j$. We label $\mathcal{B}_q^j$ so that the superscript $j$ indicates the partition $\mathcal{P}^j$ and the subscript $q$ indicates the specific cluster which is broken. Henceforth, we will call each of these vectors a \textit{breaking vector}.
The entries of $\mathcal{B}_q^j$ corresponding to nodes that do not belong to $C_q$ are set to 0, whereas those corresponding to nodes belonging to $C_q$ have the same value if they belong to the same cluster $C_i^j$. Notice that the pair $q,j$ is meaningful only if $C_q$ breaks in the transition from $\mathcal{P}$ to $\mathcal{P}^j$.
For instance, with reference to Fig. \ref{fig:balpart}, if we consider $\mathcal{P}^1$ and cluster $C_2 \in \mathcal{P}$, the breaking vector $\mathcal{B}_2^1$ has the following structure:
$\mathcal{B}_2^1 = [\underbrace{0 \quad 0}_{C_1} \quad \underbrace{\overbrace{a}^{C_2^1} \overbrace{\quad b \quad b}^{C_3^1}}_{C_2}]$.
By contrast, the breaking vector $\mathcal{B}_1^4$ is meaningless, since $C_1$ does not break in the transition from $\mathcal{P}$ to $\mathcal{P}^4$.

We define the \textit{intertwining index} $n_q^j$ as the number of clusters $C_p \in \mathcal{P}$ other than $C_q$ (i.e., $p \neq q$) which are broken after breaking cluster $C_q$ as described by vector $\mathcal{B}_q^j$. In other words, $n_q^j$ is the number of clusters of $\mathcal{P}$ (in addition to $C_q$) which are affected by the breaking action of $\mathcal{B}_q^j$, namely the number of other clusters of $\mathcal{P}$ whose stability depends on the desynchronization of $C_q$. 
For instance, in the example of Fig. \ref{fig:balpart}, $n_2^1 = 0$, because $C_1$ remains unaltered in the transition from $\mathcal{P}$ to $\mathcal{P}^1$ (notice that in $\mathcal{P}^1$ $C_1$ is labeled $C_1^1$). This means that the loss of stability (i.e., the breaking) of $C_2$ does not influence any other cluster.
By contrast, $n_1^5 = 1$, because in the transition from $\mathcal{P}$ to $\mathcal{P}^5$ or to $\mathcal{P}^6$ also $C_2$ is broken by $\mathcal{B}_1^5$. This means that the breaking of $C_1$ influences the stability of $C_2$, therefore $C_2$ is one-way dependent on $C_1$.

The breaking vectors play a fundamental role not only in checking the presence of one-way dependent clusters, but also in finding the irreducible transformation for directed networks, as described in the following.

\subsection*{Constructive method to obtain an irreducible transformation for directed networks}
As stated above, we aim to find the matrix $T$ that converts the node coordinate system to the irreducible coordinate system, thus evidencing the inter-dependencies among the perturbation components.

For undirected networks, this matrix can be found as described in \cite{siddique:2018,dellarossa:2020,zhang:2020}. 
To the best of our knowledge, there is no general method for finding this transformation matrix for directed networks; the method proposed in \cite{cho:2017} allows finding a matrix $T$, but in general it does not correspond to an irreducible transformation. Here we propose a method that outputs the irreducible transformation for directed networks and holds also (as a particular case) for undirected networks.

The key variational equation for studying cluster stability has the following structure \cite{lodi:2020}:
\begin{equation}\label{eq:variationalcompact}
\dot\deta = \Psi_1[\bs(t)] \; \deta(t) + \sum_{k=1}^{L}\Psi_2^k[\bs(t-\delta_k)] \; \deta(t-\delta_k),
\end{equation}
where $\Psi_1$ and $\Psi_2^k$ are time-dependent matrices that depend on the synchronous solution $\bs(t)$.

The block-diagonal matrix $\Psi_1$ depends only on the network (not on the coordinate change) and, through its action, $\dot{\deta}_j$ only depends on $\deta_j$. By contrast, each matrix $\Psi_2^k$ depends on the transformation matrix $T$ through a matrix $B^k$ and therefore $\dot{\deta}_j$ is related to the other perturbation components by the second term on the right-hand side of the above equation. 
Taking into account that $T$ is an orthogonal matrix (i.e., $T^{-1}=T^T$), each matrix $B^k$ has the following structure:
\begin{equation}\label{eq:Bkstruct1}
B^k = T A^k T^{-1} = \begin{bmatrix} T_\parallel \\ T_\perp
\end{bmatrix} A^k \begin{bmatrix} T_\parallel^T \quad T_\perp^T
\end{bmatrix} = \left[ \begin{array}{c;{2pt/2pt}c}
T_\parallel A^k T_\parallel^T   &  T_\parallel A^k T_\perp^T \\ \hdashline[2pt/2pt]
T_\perp A^k T_\parallel^T & T_\perp A^k T_\perp^T
\end{array}  \right ] ,
\end{equation}
where the block $T_\perp A^k T_\parallel^T$ contains only 0 entries (the proof of this statement is in the section \textit{Methods}), whereas the block $T_\parallel A^k T_\perp^T$ is null for undirected networks and contains non-zero entries for directed networks \cite{lodi:2020}.
Therefore, in general, $B^k$ can be written as
\begin{equation}\label{eq:Bkstruct}
B^k = \left[ \begin{array}{c;{2pt/2pt}c}
\multicolumn{2}{c}{B^k_\parallel} \\ \hdashline[2pt/2pt]
0 & B^k_\perp
\end{array}  \right ] .
\end{equation}
This implies that the transverse perturbations ($\deta_{Q+1},\ldots,\deta_N$) accounting for cluster stability do not depend on the parallel perturbations ($\deta_1,\ldots,\deta_Q$) along the synchronous manifold.

In particular, our main goal is to analyze the CS by decoupling the transverse perturbations \textit{as much as possible}, with the goal of reducing the dimensionality of the problem. Therefore, we look for the transformation $T$ that makes the matrix $B_\perp = \sum_{k=1}^L B^k_\perp$ block-diagonal or, at least, block-upper triangular. 
The method to construct the matrix $T$ with these characteristics is based on the breaking vectors defined above, as we explain below.

The choice of the $T$ rows related to the equitable cluster $C_q$ does not depend on the other clusters: as a matter of fact, owing to the constraints (A) and (B), the rows of $T$ related to a cluster are orthonormal to those related to any other cluster.

We construct the matrix $T$ by stacking together $Q$ matrices $T_q$ (one for each cluster), with dimensions $N_q \times N$.
The first row of matrix $T_q$ has entry $1/\sqrt{N_q}$ for each column corresponding to a node in $C_q$, zero otherwise.
The other $N_q-1$ rows are taken among the vectors $\mathcal{B}_q^j$ such that (i) the rows of $T_q$ are orthonormal and (ii) $n_q^j$ is minimum. If there is more than one vector $\mathcal{B}_q^j$ fulfilling these conditions, we select the one with the largest number of occurrences.

Once we have all matrices $T_q$ ($q=1,\ldots,Q$),
\begin{itemize}
	\item $T_\parallel$ contains the first rows of all the matrices $T_q$;
	\item $T_\perp$ contains all the other rows of $T_q$, sorted in order to make the matrix $B_\perp$ block-diagonal with minimal block sizes and with sizes decreasing along the diagonal. This is done by applying the Cuthill-McKee algorithm \cite{cuthill:1969} to $B$ and associated swaps of rows and columns to $T$.
\end{itemize}

This choice ensures that $T$ is the transformation we are looking for, because $B_\perp$ is in a form that is \emph{as block-diagonal as possible} (diagonal if $n_q^j=0$ for any $q$ and $j$).

To illustrate the method, we detail the construction of the matrix $T$ for the network of Fig. \ref{fig:balpart}.

Table \ref{tab:esempio_rotture_intertwIndex} shows the breaking vectors $\mathcal{B}_q^j$ and the corresponding intertwining indices $n_q^j$ for each partition $\mathcal{P}^j$ and for the equitable clusters $C_1$ and $C_2$ of partition $\mathcal{P}$ in Fig. \ref{fig:balpart}.

\begin{table}[h!]
	\centering
	\begin{tabular}{ | c | c| } 
		$\mathcal{B}_1^5 = [\overbrace{\underbrace{a}_{C^5_1} \quad \underbrace{b}_{C^5_4}}^{C_1} \quad \overbrace{0 \quad 0 \quad 0}^{C_2}]$ & $n_1^5 = 1$ \\  \hline \\
		$\mathcal{B}_1^6 =  [\underbrace{a}_{C^6_1} \quad \underbrace{b}_{C^6_2} \quad 0 \quad 0 \quad 0]$ & $n_1^6 = 1$ \\
		\hline \hline \\
		$\mathcal{B}_2^1 = [0 \quad 0 \quad \underbrace{a}_{C^1_3} \quad \underbrace{b \quad b}_{C^1_2}]$ & $n_2^1 = 0$\\ \hline \\
		$\mathcal{B}_2^2 = [0 \quad 0 \quad \underbrace{a \quad a}_{C^2_2} \quad \underbrace{b}_{C^2_3}]$ & $n_2^2 = 0$\\ \hline
		$\mathcal{B}_2^3 = [0 \quad 0 \quad a \quad b \quad a]$ & $n_2^3 = 0$\\ \hline
		$\mathcal{B}_2^4 = [0 \quad 0 \quad a \quad b \quad c]$ & $n_2^4 = 0$\\ \hline
		$\mathcal{B}_2^5 = [0 \quad 0 \quad a \quad b \quad a]$ & $n_2^5 = 1$ \\  \hline
		$\mathcal{B}_2^6 = [0 \quad 0 \quad a \quad b \quad c]$ & $n_2^6 = 1$\\ \hline
	\end{tabular}
	\caption{\footnotesize{Breaking vectors and intertwining indices. Breaking vectors $\mathcal{B}_q^j$ and corresponding intertwining indices $n_q^j$ for each partition $\mathcal{P}^j$ and for the equitable clusters $C_1$ and $C_2$ of the network of Fig. \ref{fig:balpart}, partition $\mathcal{P}$. For the first four breaking vectors also the broken clusters are shown. Cluster $C_1$ is broken only in partitions $\mathcal{P}^5$ and $\mathcal{P}^6$ (see the first two rows and Fig. \ref{fig:balpart}). In $\mathcal{P}^5$ (see Fig. \ref{fig:balpart}), in particular, $C_1$ is broken as described in the first row ($\mathcal{B}_1^5$) and also $C_2$ is broken, therefore $n_1^5 = 1$. Cluster $C_2$ is broken in all partitions. In particular,  in $\mathcal{P}^1$ it is broken as described in the third row ($\mathcal{B}_2^1$) and $C_1$ is not broken (see Fig. \ref{fig:balpart}), therefore $n_2^1 = 0$.}}
	\label{tab:esempio_rotture_intertwIndex}
\end{table}

The matrices $T_q$ ($q=1,2$) are built as follows.
Here we detail the construction of the matrix $T_2$, by taking the following steps:
\begin{enumerate}
	\item We initialize $T_2$ with the row corresponding to the synchronous manifold: $T_2 = 1/\sqrt{3} \begin{bmatrix}
	0 & 0 & 1 & 1 & 1
	\end{bmatrix}$.
	\item Since there is more than one breaking vector fulfilling conditions (i) and (ii), among all the breaking vectors $\mathcal{B}_2^j$ with $n_2^j = 0$, we select and append to $T_2$ the ones ($\mathcal{B}_2^3$ and $\mathcal{B}_2^4$) whose structure appears more frequently (twice) in Tab. \ref{tab:esempio_rotture_intertwIndex} and that are orthonormal to the other rows of $T_2$:
	\begin{equation*}
	T_2 = \begin{bmatrix}
	0 & 0 & \frac{1}{\sqrt{3}} & \frac{1}{\sqrt{3}} & \frac{1}{\sqrt{3}} \\
	0 & 0 & -\frac{1}{\sqrt{6}} & \frac{2}{\sqrt{6}} & -\frac{1}{\sqrt{6}} \\
	0 & 0 & -\frac{1}{\sqrt{2}} & 0 & \frac{1}{\sqrt{2}} \\
	\end{bmatrix}    
	\end{equation*}
	\item Since $T_2$ has already the proper size ($N_2 = 3$ rows), we stop here. In case this condition was not satisfied, we would need to consider the breaking vectors with larger intertwining indices $n_2^j$.
\end{enumerate}

$T_1$ is obtained analogously and is as follows:
\begin{equation*}
T_1 = \frac{1}{\sqrt{2}}\begin{bmatrix}
1 & 1 & 0 & 0 & 0 \\
1 & -1 & 0 & 0 & 0 
\end{bmatrix}
\end{equation*}

Now we can build up the matrix $T$: $T_\parallel$ contains the first row of $T_1$ and the first row of $T_2$, $T_\perp$ contains all the other rows of $T_1$ and $T_2$. The rows of $T_\perp$ are sorted in order to make the matrix $B=TA^1T^T$ as close as possible to a diagonal matrix:
\begin{equation}
T = \begin{bmatrix} T_\parallel \\ \hdashline[3pt/4pt] T_\perp
\end{bmatrix} = \begin{bmatrix}
\frac{1}{\sqrt{2}} & \frac{1}{\sqrt{2}} & 0	 & 0 & 	0\\
0 & 0 & \frac{1}{\sqrt{3}} & \frac{1}{\sqrt{3}} & \frac{1}{\sqrt{3}}\\ \hdashline[3pt/4pt]
0 & 0 & -\frac{1}{\sqrt{6}} & \frac{2}{\sqrt{6}} & -\frac{1}{\sqrt{6}}\\
-\frac{1}{\sqrt{2}} & \frac{1}{\sqrt{2}} & 0 & 0 & 0\\
0 & 0 & -\frac{1}{\sqrt{2}} & 0 & \frac{1}{\sqrt{2}}
\end{bmatrix}
\begin{array}{l}
\rightarrow C_1\\
\rightarrow C_2\\
\rightarrow C_2\\
\rightarrow C_1\\
\rightarrow C_2
\end{array}.
\end{equation}

\begin{figure}[h]
	\centering
	\includegraphics[width = 0.45\linewidth]{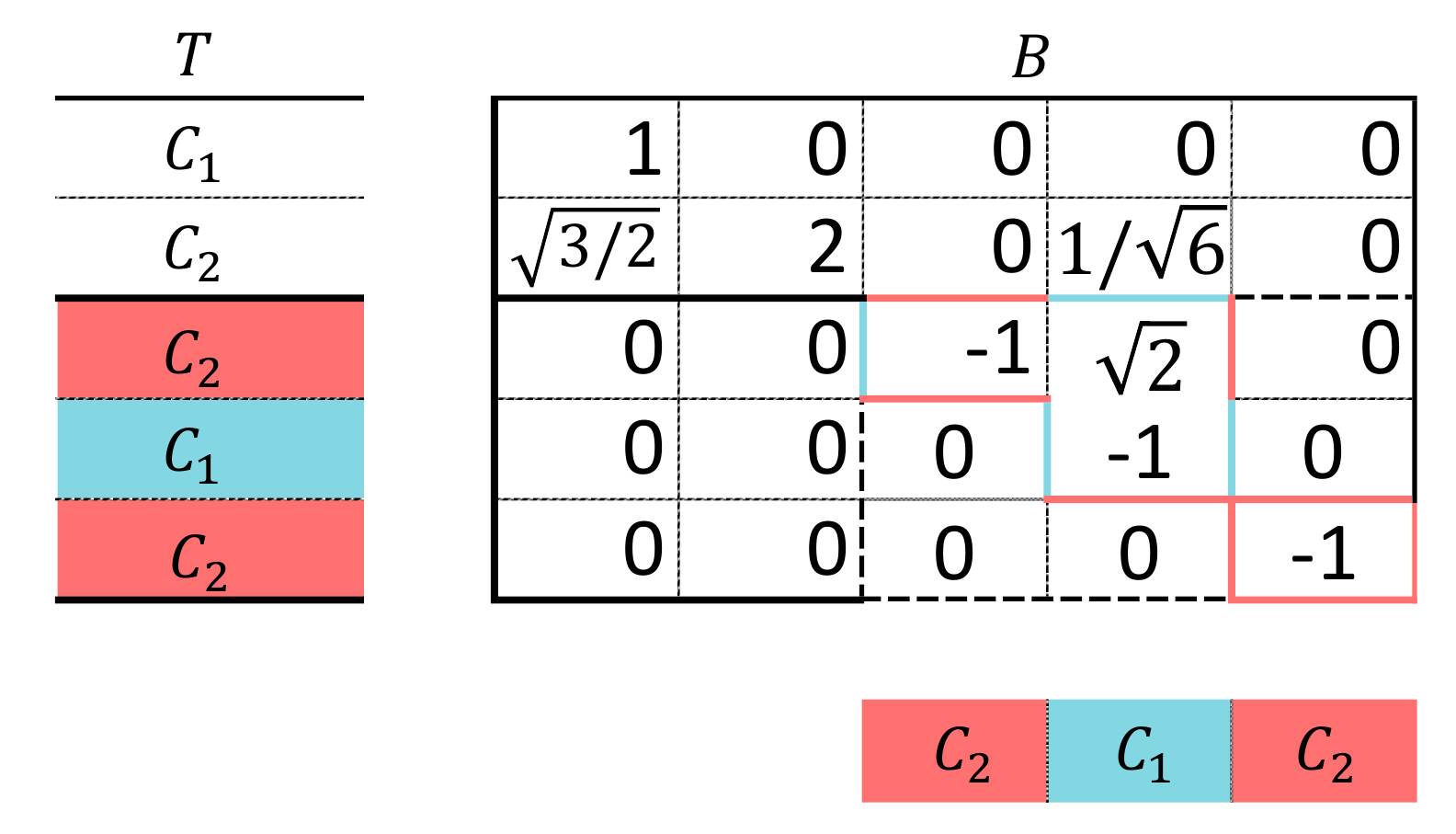}
	\caption{\footnotesize{Matrices $T$ and $B$ for the example of Fig. \ref{fig:balpart}, partition $\mathcal{P}$. Structure of the matrix $T$ (left, with the rows of $T_\perp$ colored as the corresponding cluster) and complete matrix $B$ (right). The transverse block $B_\perp$ (within the dashed rectangle) has an upper triangular structure, containing two sub-blocks with non-zero entries (color-bordered blocks). The upper sub-block (with red and blue borders) characterizes (see the corresponding rows in the $T$ structure) the perturbation affecting clusters $C_2$ (first row of the sub-block), which depends on both clusters, and the one affecting $C_1$ (second row), which depends only on $C_1$ (single diagonal entry). The red rectangle evidences a second perturbation affecting cluster $C_2$, which is independent of cluster $C_1$. Therefore, $C_1$ is an independent cluster, whereas $C_2$ is one-way dependent on $C_1$.}}
	\label{fig:matrice_B_esempio}
\end{figure}

The corresponding matrix $B$ is shown in Fig. \ref{fig:matrice_B_esempio}.
From the block-upper triangular structure of the transverse block $B_\perp$ in Fig. \ref{fig:matrice_B_esempio}, we see that the perturbation affecting cluster $C_1$ is independent of those affecting cluster $C_2$. One of the two perturbations affecting cluster $C_2$ is in turn independent of that affecting cluster $C_1$. Assuming that this perturbation is stable, the remaining perturbation affecting cluster $C_2$ is also independent. This reflects the fact that the stability of $C_1$ does not depend on the stability of $C_2$, whereas the stability of $C_2$ depends on the stability of $C_1$. Therefore, as stated above, $C_2$ is one-way dependent on $C_1$.

In summary (see Fig. 1 in the \textit{Supplementary Material}), 
the characteristics of inter-dependence between clusters (intertwined, one-way dependent or independent) in the general class of networks \eqref{eq:dyneq} can be deduced from the ISS and from the breaking vectors; therefore, they are based exclusively on the network topology. Through the matrix $B$ related to the irreducible transformation matrix $T$ built up as described above, we can evaluate the stability of clusters, thus extending the method proposed in \cite{lodi:2020} to a larger class of networks. In so doing, we solve an open problem \cite{zhang:2020}.
Moreover, by using the irreducible transformation matrix $T$ based on the breaking vectors, one can easily detect the presence of intertwined or one-way dependent clusters by inspecting the submatrix $B_\perp$:
\begin{itemize}
	\item if $B_\perp$ is block-upper triangular, there is a one-way dependency of the clusters involved in a specific block;
	\item otherwise, if $B_\perp$ is block-diagonal, the clusters involved in a specific block are intertwined.
\end{itemize}

\begin{figure}[h]
	\centering
	\includegraphics[width = 1\linewidth]{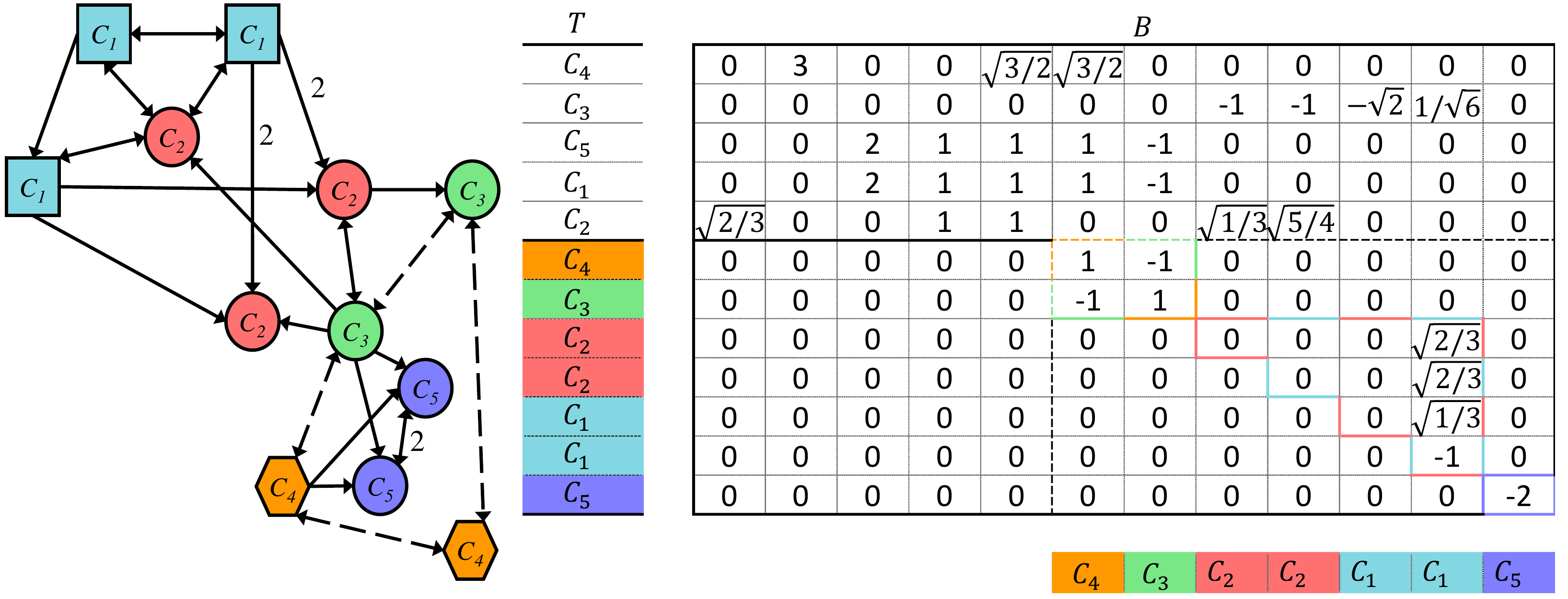}
	\caption{\footnotesize{Summarizing illustrative example. The minimal balanced partition (left panel) for the network of Fig. \ref{fig:exnetw}A: it has $Q=5$ clusters $C_1 = \{1,2,3 \}$ (colored blue), $C_2 = \{4,5,6 \}$ (colored red), $C_3 = \{7,8 \}$ (colored green), $C_4 = \{9,10 \}$ (colored orange), and $C_5 = \{11,12 \}$ (colored purple). Corresponding matrices (right panel): structure of $T$ (left) and complete matrix $B$ (right). The $T$ rows corresponding to $T_\perp$ are colored as the corresponding cluster. The transverse block $B_\perp$ (within the dashed rectangle) has an upper triangular structure, containing three sub-blocks with non-zero entries (color-bordered blocks). The purple rectangle characterizes the only perturbation affecting cluster $C_5$ (see the corresponding row in the $T$ structure), which is independent (single diagonal entry). The green-orange $2 \times 2$ rectangle evidence the two perturbations affecting clusters $C_3$ (second row of the block), which depends also on cluster $C_4$, and $C_4$ (first row of the block), which depends also on cluster $C_3$; this means that clusters $C_3$ and $C_4$ are intertwined. The blue-red upper triangular block is concerned with the two perturbations affecting clusters $C_1$ and $C_2$: both of them depend only on $C_1$, meaning that $C_2$ is one-way dependent on $C_1$.}}
	\label{fig:excompl}
\end{figure}

To better appreciate the advantages of detecting the inter-dependencies among clusters in more complex networks by inspecting $B_\perp$, Fig. \ref{fig:excompl} shows the application of the method to the network of Fig. \ref{fig:exnetw}A, whose minimal balanced partition has $Q=5$ clusters.

In the following, the proposed method is applied to two examples.

\subsection*{Violin players}
The first example is from \cite{shahal:2020} and is concerned with a human network of violin players executing a musical phrase. The coupling strength between two players is defined as the ratio between the sound volume of the coupled violin and the sound volume of the player's own violin, while maintaining constant the total volume that each player hears. The dynamics of the phase $\varphi_i$ of the $i$-th violin player is simulated according to the model proposed in \cite{shahal:2020}:
\begin{equation}
\dot{\varphi}_i = \omega_i + \kappa \sum_{j=1}^{N} A_{ij} \sin \left( \varphi_j (t-\delta) - \varphi_i (t) \right) ,
\end{equation}
where $\omega_i$ is the eigenfrequency of the $i$-th player, $\kappa = 0.2$ Hz is the coupling strength and $\delta$ is the delay in seconds between the players. Each player can decide to ignore one input, in which case a directed connection is established. In particular, in \cite{shahal:2020} for some delay values each player spontaneously decides to ignore one input, leading to the formation in the described experiments of directed networks of particular relevance.

Unlike \cite{shahal:2020}, where the eigenfrequencies were randomly distributed in the range $[0.25,0.3]$ Hz, here we consider all players ($i=1,\ldots,N$) with the same eigenfrequency $\omega_i = 0.25$ Hz.
We focus on some of the ring topologies with $N=8$ violin players studied in \cite{shahal:2020}. On the basis of the delay $\delta$, the experiments in \cite{shahal:2020} have determined that each player either pays attention to both closest neighbors or decides to ignore either one of them, which corresponds to effectively changing the network connectivity.

\begin{figure}[h!]
	\centering
	\includegraphics[width=\linewidth]{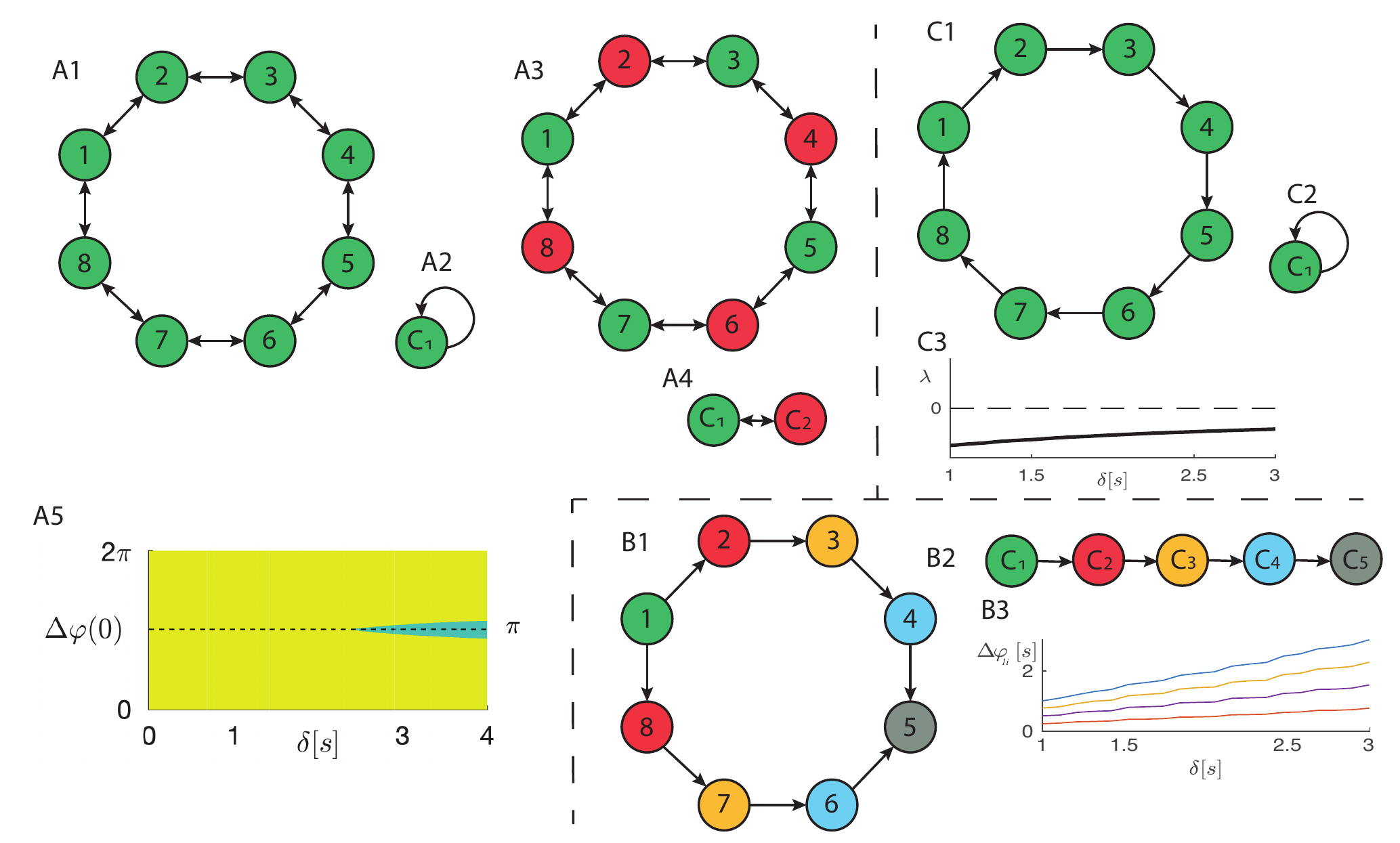}
	\label{fig:outofphase_8}
	\caption{\footnotesize{Example 1: ring of 8 violin players. For delay $\delta < 1$s or $3 < \delta < 4$s, the network of violin players can be represented as undirected (panels A), and for 1s $< \delta < 3$s as a directed network with `arrowhead' topology, oriented either clockwise (as in panels B and C) or counterclockwise. Undirected network with minimal balanced coloring (A1, `in-phase' global synchronization pattern) and with (non-minimal) balanced coloring (A3) and corresponding quotient networks (A2 and A4, respectively); basins of attraction (A5) of the two stable solutions admitted by the network A4: `in-phase' global synchronization (yellow region) and `out-of-phase' synchronization with two clusters oscillating in anti-phase (green region). Directed network with `arrowhead' topology and (non-minimal) balanced coloring (B1), corresponding quotient network (B2), and asymptotic phase differences $\Delta \varphi_{i1}$ vs $\delta$ (B3): it is apparent that $\Delta \varphi_{i1}$ grows linearly with $\delta$. Directed network with `arrowhead' topology and minimal balanced coloring (C1) and corresponding quotient network (C2); MLE $\lambda$ corresponding to this synchronization pattern (C3), showing stability as $\lambda(\delta) < 0$ for any $\delta \in (1,3)$.}} 
\end{figure}

\textit{a) Undirected network} -- With low delay (roughly, $\delta < 1$s), each player hears both closest neighbors, which corresponds to having bidirectional connections, as shown in panel A1. The corresponding minimal balanced coloring is uniform, meaning that the quotient network contains only one node (panel A2) and the network can synchronize globally. This synchronization pattern corresponds to the eight violins being played `in-phase' and is stable because the corresponding MLE is negative for any value of $\delta$.

Similarly, in the case of large delay $\delta \in [3,4]$s, each player listens to both closest neighbors, thus leading again to global synchronization. However, an open question is whether the network in panel A1 can admit the `out-of-phase' synchronization  experimentally found in \cite{shahal:2020}. To this end, we consider the balanced coloring shown in panel A3, corresponding to two clusters and to a quotient network of two nodes (panel A4). Let $\Delta \varphi$ be the phase lag between the two clusters. For each $\delta$ value, the quotient network has two stable equilibrium points: one for $\Delta \varphi = 0$ and one for $\Delta \varphi = \pi$. The first one corresponds to the global synchronization predicted by the minimal balanced coloring, whereas the second one corresponds to a stable out-of-phase synchronization with two clusters oscillating in anti-phase. Panel A5 shows the pairs $(\delta, \Delta \varphi(0))$ that generate trajectories eventually approaching either one of the two equilibrium points: $\Delta \varphi = 0$ (yellow points) and $\Delta \varphi = \pi$ (green points). The figure evidences the presence of multistability, which is in agreement with the analysis in \cite{shahal:2020}, even if anti-phase synchronization was not found with the original model. As a matter of fact, in \cite{shahal:2020} $\delta$ is increased linearly with time and the initial condition $\Delta(0)$ is therefore always close to 0, thus preventing the network from reaching the experimentally found stable anti-phase pattern. There are some quantitative discrepancies between Ref. \cite{shahal:2020} and our work, regarding the delay values at which the out-of-phase state of synchronization starts to be experimentally detected (less than 2s). The quantitative differences between experimental and numerical results are mainly due to the simple model adopted in both papers. The further quantitative differences in the simulated results can be ascribed to two main changes in our modelling approach with respect to \cite{shahal:2020}:  (i) we consider the asymptotic network state for a fixed delay, whereas \cite{shahal:2020} considers a time-varying delay, namely the results contain also transient information; (ii) we consider the same oscillation frequency for all the oscillators, in order to allow for the existence of cluster synchronization. We also remark that detection of multistability is harder in the original network, because the space of the phase lags is $\mathbb{R}^7$, whereas in the quotient network it is only $\mathbb{R}^1$, owing to the dimensionality reduction. Notice that for low $\delta$ values the equilibrium point $\Delta \varphi = \pi$ is unstable, indicating that we can only have global synchronization of the violin players. 

The matrix $B$ for this network is reported and analyzed in the \textit{Supplementary Material}, Sec. 3.1.

\textit{b) Directed network} -- If $\delta \in [1,3]$s, each violin player ignores one of its closest neighbors, leading to either one of two different directed networks.
The first network is configured according to the `arrowhead' topology shown in panel B1. The corresponding minimal balanced coloring (panel B2) contains 5 colors, therefore the network admits $Q=5$ clusters (see the quotient network in panel B2), two of which are trivial, which excludes the possibility of global synchronization. The matrix $B$ for this network is reported and analyzed in the \textit{Supplementary Material}, Sec. 3.2. The corresponding clusters are stable. The phase difference $\Delta\varphi_{ij} = \varphi_i - \varphi_j$ between two clusters grows linearly with $\delta$, as shown in panel B3 using node 1 as reference (i.e., $j=1$).

The second network has unidirectional connections with the same orientation (either clockwise or counterclockwise) which corresponds to the uniform balanced coloring in panel C1 and global synchronization associated with a rotational symmetry. Here we do not obtain the `vortex' pattern found in \cite{shahal:2020}, due to our assumption of identical eigenfrequencies $\omega_i$ for all players.

\subsection*{Neural activity and chimera states}
The second example is adapted from \cite{majhi:2017}.
\begin{figure}
	\centering
	\includegraphics[width=0.75\linewidth]{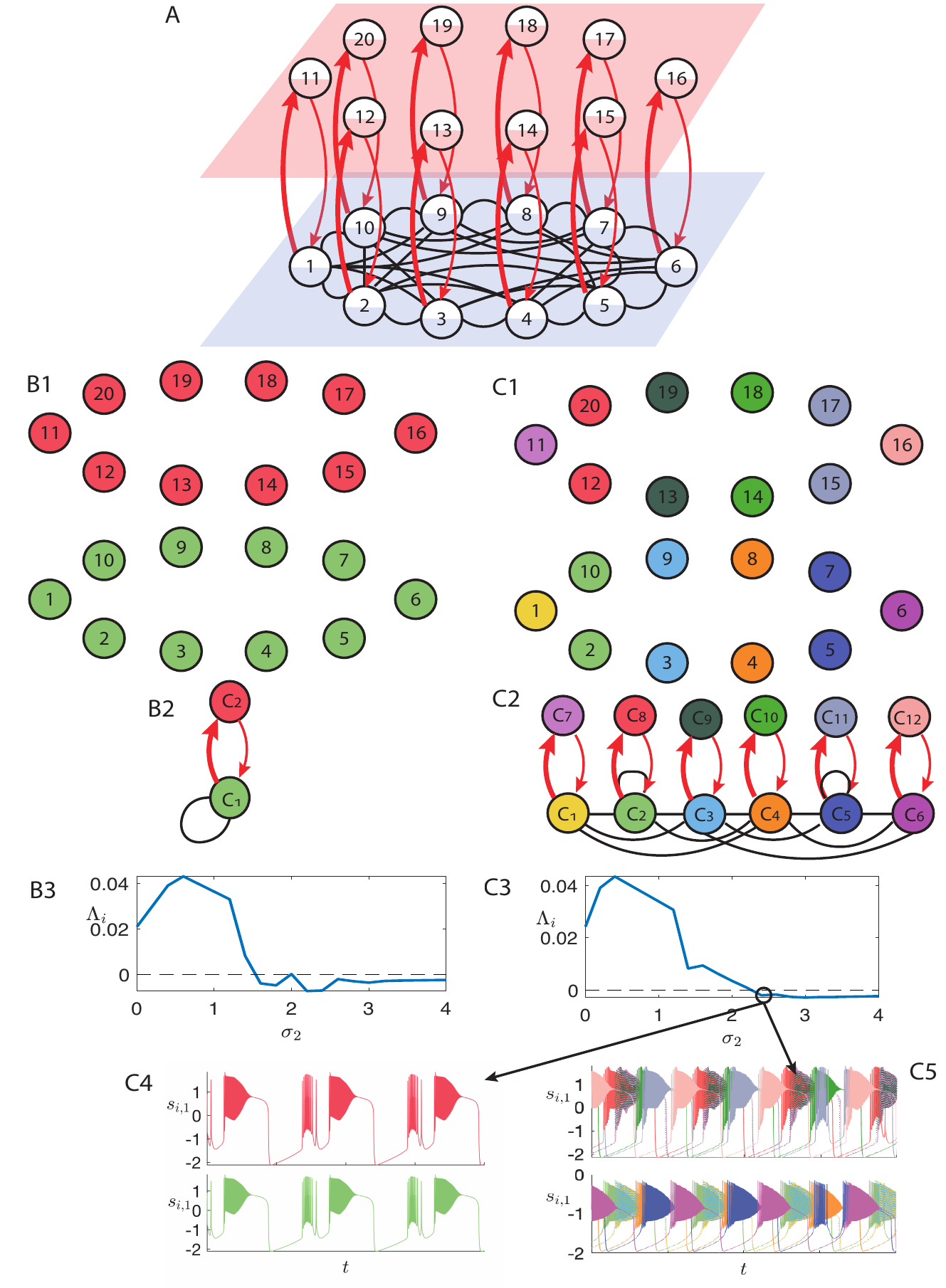}
	\caption{\footnotesize{Example adapted from \cite{majhi:2017}. (A) Multilayer network with $N=20$ nodes, $L=2$ kinds of links, undirected intra-layer and directed layer-to-layer connections. (B1) Minimal balanced coloring/partition of the network, corresponding to $Q=2$ clusters $C_1$ (green nodes, from 1 to 10) and $C_2$ (red nodes, from 11 to 20). (B2) Corresponding quotient network. (B3) MLEs $\Lambda_q$ corresponding to cluster $C_q$ ($q=1,2$) vs. $\sigma_2$. (C1) Non-minimal balanced coloring/partition, corresponding to $Q=12$ clusters, four of which are trivial (nodes 1, 6, 11, 16), while the other eight contain two nodes (2-10, 3-9, and so forth). (C2) Corresponding quotient network. (C3) MLEs $\Lambda_q$ corresponding to cluster $C_q$ ($q=1,\ldots,Q$) vs. $\sigma_2$. (C4,C5) Time plots (first component of the synchronous cluster states $\bs_{i}$ for $i=1,\ldots,Q$) evidencing the presence of multistability for $\sigma_2 = 2.3$: panel C4 corresponds to the minimal balanced partition of panel B1 (two distinct traces), while panel C5 corresponds to the non-minimal balanced partition of panel C1 (twelve distinct traces).}}
	\label{fig:chimeranetw}
\end{figure}
It is concerned with a network (shown in Fig. \ref{fig:chimeranetw}, panel A) consisting of $N=20$ neurons, distributed in $M=2$ layers with the same number $N/2$ of nodes and connected through $L=2$ kinds of synapses. The upper layer (layer I, in red) is composed of uncoupled identical neurons (see section \textit{Methods}) and the lower layer (layer II, in blue) is composed of coupled identical neurons, different from those of layer I. Let $1,\dots,N/2$ be the labels of neurons of layer II and $N/2 + 1,\dots,N$ those of neurons in layer I.
Each neuron $i=1,\dots,N/2$ of layer I is connected to neuron $i+N/2$ of layer II through a double directed chemical excitatory synapse with strength $\sigma_2$ in one direction (see thick red connections from layer II to layer I) and $\sigma_2/4$ in the opposite direction (see thin red connections from layer I to layer II).

In layer II, each neuron is connected to its $P=6$ closest neighbors through identical electrical synapses (gap junctions) with strength $\sigma_1$. Nodes of layer I interact only through the layer II neurons.

The minimal balanced coloring of the network (panel B1) evidences the presence of two clusters $C_1 = \{1,\dots,N/2\}$ and  $C_2 = \{N/2+1,\dots,N\}$, namely cluster $C_1$ contains all the neurons of layer I and cluster $C_2$ all those of layer II.
The corresponding quotient network (panel B2) is composed of two nodes only. Through the proposed method (see the \textit{Supplementary Material}, Sec. 3.3, for details), we obtain the matrix $T$ corresponding to the irreducible transformation of coordinates. The corresponding matrix $B$ shows that cluster $C_2$ is intertwined with $C_1$.

Following \cite{majhi:2017,ruzzene:2020}, we study the cluster stability by changing the synaptic efficacy $\sigma_2$, with $\sigma_1 = 0.005$. Panel B3 shows the MLEs $\Lambda_q$ corresponding to cluster $C_q$ ($q=1,2$) vs. $\sigma_2$ (overlapping curves, because the clusters are intertwined): it is evident that these clusters are stable only for $\sigma_2$ values larger than a threshold approximately equal to 1.5.

Next we want to check the presence of chimera states \cite{kuramoto:2002,abrams:2004}, which are characterized by the coexistence of coherent and incoherent dynamics.
This kind of state has been detected in several neural systems and it has been associated to either pathological or normal activity \cite{majhi:2019,wang:2020}. For instance, during an epileptic seizure, some regions of the brain are strongly synchronized, whereas others are desynchronized \cite{ayala:1973}. Similarly, the two brain hemispheres of a sleeping dolphin have independently synchronized and desynchronized behaviors at the same time \cite{mukhametov:1977}: one hemisphere is in sleep and another remains awake. This \textit{unihemispheric slow-wave sleep} has also been found in other aquatic animals and migrating birds \cite{rattenborg:2000}.

The stable clusters $C_1$ and $C_2$ of panel B2 correspond to coherent dynamics. 
In order to check whether the considered network can generate other patterns of synchronization, we consider a balanced coloring that is not minimal.

The network admits many non-minimal balanced colorings, one of which is shown in panel C1, containing a larger number of clusters ($Q=12$), as pointed out by the corresponding quotient network (panel C2).
Also in this case, through the proposed method (see the \textit{Supplementary Material}, Sec. 3.3, for details), we obtain the matrix $T$ corresponding to the irreducible transformation of coordinates. Analysis of the corresponding matrix $B$ allows to conclude that all nontrivial clusters are intertwined, meaning that the breaking of one cluster implies that all other clusters also break. The possibility of observing the solution shown in either panel B1 or C1 depends on their corresponding stability, which can change with $\sigma_2$. If more than one cluster is stable for a given value of $\sigma_2$, the observed solution depends on the initial condition.

Panel C3 shows the MLEs $\Lambda_q$ corresponding to cluster $C_q$ ($q=1,\ldots,Q$) vs. $\sigma_2$ (again, all the curves overlap because the clusters are intertwined and the same MLE determines the stability of all clusters for each value of $\sigma_2$): it is evident that all these clusters are stable for $\sigma_2$ values exceeding a threshold $\hat{\sigma} \approx 2.2$.

Panels C4 and C5 (obtained for $\sigma_2 = 2.3$ but with different initial conditions) show time plots of the first component of the $q$-th cluster state $\bs_q(t)$, for layer I (upper panels) and layer II (lower panels). Traces are colored as the corresponding clusters in panels B2 and C2, respectively.
In panel C5 we see that almost each cluster corresponds to a different trace, with the exception of the pairs $\{C_1,C_3\}$ and $\{C_7,C_9\}$, which are synchronized: this means that in the original network nodes 1, 3 and 9 are in a cluster and 11, 13, 19 in a different one.
By contrast, in panel C4 the traces correspond to the minimal balanced coloring. This indicates the presence of multistability for $\sigma_2 = 2.3$. We remark that also other stable patterns may emerge: for instance, by setting initial conditions such that the states of C4 and C6 are close, also C4 and C6 synchronize.

Now we consider the same network in Fig. \ref{fig:chimeranetw}A, but we remove all interlayer connections from the upper layer to the lower layer. In this case, the matrices $T$ and $B$ do not depend on $\sigma_2$.
In both cases, we focus on the minimal balanced coloring only.
By removing the connections from the upper to the lower layer, the $2 \times 2$ submatrix $B_\perp$ is block-upper triangular and the structure of the corresponding matrix $T$ evidences that cluster $C_1$ is one-way dependent on $C_2$ (see the \textit{Supplementary Material}, Sec. 3.4, for details).
Conversely, if we remove the connections from the lower to the upper layer, i.e., the thick red connections from layer II to layer I, $B_\perp$ is again a unique $2 \times 2$ block-upper triangular submatrix and the structure of the corresponding matrix $T$ evidences that cluster $C_2$ is one-way dependent on $C_1$.

\section*{Discussion}\label{sec:discussion}
There are two main contributions of this paper. We defined one-way dependent clusters and provided a method to reduce the CS stability problem for directed networks of coupled nonlinear dynamical systems in an irreducible form.

The concept of one-way dependent clusters is different from the previously reported concept of intertwined clusters, for which the desynchronization of a cluster A is subordinate to that of another cluster B and \textit{vice versa}. As most real networks are directed, this concept has profound implications on our understanding of clustering in real-world settings. We illustrated this statement via two examples: a human network of violin players executing a musical piece, for which directed interactions may be either activated or deactivated by the musicians, and a multilayer neural network with directed layer-to-layer connections.

Our work supersedes some of the limitations of previous research in the area of cluster synchronization. Reference \cite{pecora:2014} proposed a method that was effective in decoupling the CS stability problem in a number of irreducible lower-dimensional problems for arbitrary networks. However, the method only worked for \textit{orbital clusters} and not also for \textit{equitable clusters}. It is known \cite{bari:2006} that all orbital clusters are necessarily equitable but the converse is not necessarily true, so there exist equitable cluster partitions that are not orbital.  There are notable examples of networks for which the minimal orbital cluster partition and the minimal equitable cluster partition do not coincide 	\cite{kudose:2009}. Reference \cite{cho:2017} focused on undirected and unweighted networks, including the case of equitable clusters, with the possibility of extending this theory to the case of directed and weighted networks. However, this method is unable to guarantee a decomposition of the stability problem in irreducible blocks. In \cite{irving:2012,zhang:2020} an algebra-based approach is adopted to analyze the stability of synchronized clusters in undirected networks, as an alternative to the methods based on the theory of groups or groupoids \cite{golubitsky:2006,pecora:2014}. In particular, the method presented in \cite{zhang:2020} is able to output irreducible blocks defined in a mathematically meaningful way, but does not apply to the case of directed networks. Explicit analytical conditions for the (local) stability of synchronized clusters in sparse and weighted networks of heterogeneous Kuramoto oscillators have been proposed in Ref. \cite{menara:2020}. To conclude, to the best of our knowledge, our method is the first one that permits to reduce the stability problem in a number of irreducible lower-dimensional subproblems for the important class of directed networks. As the model for the node dynamics is quite general \cite{barzel:2013} and also the topology is now general (directed networks), the proposed method can be applied to a broad class of dynamical networks, overcoming the current limitations pointed out above. Moreover, the proposed examples show the high flexibility of our approach.

Our method to build up the irreducible transformation matrix $T$ and the corresponding matrix $B$ is based in part on the methods for the stability analysis of CS solutions \cite{pecora:2014,sorrentino:2016,schaub:2016,siddique:2018,dellarossa:2020,lodi:2020} and in part on the methods for finding and analyzing ISS in networks of dynamical systems \cite{stewart:2003,aguiar:2009,aguiar:2011,aguiar:2014,aguiar:2017,aguiar:2018,neuberger:2020}. In reference to the latter class of methods, in particular, the correspondence between cluster synchronization and the concepts of ISS and balanced partition was already established in \cite{aguiar:2009,aguiar:2018}. In \cite{aguiar:2018} synchrony-breaking codimension-one bifurcations were analyzed using the ISS in weighted regular networks, namely networks in which the sum of the weights of the edges directed to the node is constant.
However, these methods have been applied to classes of networks less general than in this paper, even if reference \cite{aguiar:2014} (which focuses on regular networks) provides some clues about generalization to non-regular and non-homogeneous networks. We remark that the use of these tools for the stability analysis of the general class of networks \eqref{eq:dyneq} with optimal separation of the perturbation modes and for the associated construction of matrices $T$ and $B$ is new.

One of the main limitations of the proposed approach is its computational complexity, mainly due to the algorithm for finding the ISS \cite{kamei:2013,aguiar:2014,neuberger:2020}.
As stated above, in the case of undirected networks, some methods (based either on the irreducible representation of the symmetry group or on simultaneous block diagonalization of matrices) have been successfully applied to reduce the computational complexity of computing the transformation matrix T. Reducing the computational burden in the case of directed networks remains an open challenge.

Our method provides a key for characterizing model-based networks of dynamical systems and for analyzing CS (for networks of persistently oscillating elements, as in the case studies) or \textit{consensus} \cite{sorrentino:2020group} (if the nodes of the network converge to an equilibrium point instead of a periodic solution) in directed networks. Among the many possible applications cited in the introduction, we believe that one of the most promising is in the area of neuroscience. The field of network neuroscience \cite{bassett:2017} is one of the emerging frontiers in the investigation of complex brain functions. Principles and mechanisms underlying some of these functions may be related to the emergence of cluster synchronization in directed networks.

\section*{Methods}

\subsection*{Invariant synchrony subspaces}
A subspace $U$ of $\mathbb{R}^n$ is invariant under a matrix $A \in \mathbb{R}^{n \times n}$ if $AU = \{Au: u \in U\}$ is a subset of $U$.

The ISS have been found following the theory developed in \cite{stewart:2003,aguiar:2011,aguiar:2014,aguiar:2017,aguiar:2018} and using the algorithm proposed in \cite{neuberger:2020}.

\subsection*{Proof of Eq. \ref{eq:Bkstruct}}
We first recall that the cluster synchronization manifold is $A^k$-invariant \cite{schaub:2016}: this means that the image through $A^k$ of any vector that lies in the cluster-synchronous subspace (any linear combination of the first $Q$ rows of the matrix $T$) lies in the space spanned by the first $Q$ rows of the matrix $T$. Since
\[
B^k=T A^k T^{-1}
\]
and $T$ is orthonormal, we obtain
\[
A^k T^T = T^T B^k = \left((B^k)^T T\right)^T.
\]

Recalling that the remaining $N-Q$ rows of $T$ are orthogonal to the first $Q$, $(B^k)^T$ must have 0 entries in all the last $N-Q$ columns of the first $Q$ rows. This means that the block $T_\perp A^k T_\parallel^T$ in Eq. \eqref{eq:Bkstruct1} is null and therefore $B^k$ is block-upper triangular.

\subsection*{Models of example 2}
Largely following \cite{majhi:2017}, the equations governing the dynamics of layer I are
\begin{equation}
\begin{array}{l}
\dot{x}_{\frac{N}{2}+i} = a x_{\frac{N}{2}+i}^2 - x_{\frac{N}{2}+i}^3 - y_{\frac{N}{2}+i} - z_{\frac{N}{2}+i} + \sigma_2 (d - x_{\frac{N}{2}+i}) h^2(x_{i})\\
\dot{y}_{\frac{N}{2}+i} = (a + \alpha) x_{\frac{N}{2}+i}^2 - y_{\frac{N}{2}+i}\\
\dot{z}_{\frac{N}{2}+i} = c(b x_{\frac{N}{2}+i} - z_{\frac{N}{2}+i} + e)
\end{array},
\end{equation}
where $a = 2.8$, $\alpha = 1.6$, $b = 9$, $c = 0.001$, $d = 2$, $e = 5$,
\begin{equation*}
h^2(x) = \frac{1}{1+e^{-\lambda (x - \Theta)}}
,
\end{equation*}
$\lambda = 10$, $\Theta = -0.25$.

Similarly, the equations governing the dynamics of layer II are the following:
\begin{equation}
\begin{array}{l}
\dot{x}_{i} = a x_{i}^2 - x_{i}^3 - y_{i} - z_{i} + \sigma_1 \sum_{j \in N_S} h^1(x_{i},x_{j}) + \sigma_2 (d - x_{i}) h^2(x_{\frac{N}{2}+i})\\
\dot{y}_{i} = (a + \alpha) x_{i}^2 - y_{i}\\
\dot{z}_{i} = c(b x_{i} - z_{i} + e)
\end{array},
\end{equation}
where $h^1(x_{i},x_{j}) = x_j - x_{i}$ and $N_S$ is the neighborhood set of the $P$ neurons closest to the $i$-th one and the other parameters are set to the same values as above, with the exception of $\alpha = 1.7$, in order to have different neuron models in the two layers.

\section*{Data availability}
All data generated or analysed during this study are included in this published article (and its supplementary information files).

\section*{Code availability}
The source code for the numerical simulations presented in the paper will be made available upon request.


\section*{Acknowledgments}
The authors would like to express their sincere appreciation to Fabio Della Rossa for many useful inputs and valuable comments. They also thank Martin Golubitsky, Mauro Parodi and Louis Pecora for reading the first version of this manuscript and providing positive and constructive comments.

\section*{Author contributions}
F.S. and M.S. designed and supervised the research; M.L. performed the research; M.L., F.S. and M.S. analyzed the data and interpreted the results; M.S. wrote and revised the manuscript; M.L. and F.S. contributed to writing and revised the manuscript.

\section*{Competing interests}
The authors declare no competing interests.

\end{document}